\documentclass{ps}
\usepackage{amsfonts}
\usepackage{graphicx}
\usepackage{amsmath}
\usepackage{amssymb}%
\providecommand{\U}[1]{\protect\rule{.1in}{.1in}}
\newtheorem {theorem}{Theorem}[section]
\newtheorem {proposition}{Proposition}[section]
\newtheorem {corollary}{Corollary}[section]

\newtheorem{lemma}{Lemma}[section]

\begin{document}
\title{Further Refinement of Self-normalized Cram\a'{e}r-type moderate Deviations}
\author{Hailin Sang}\address{Department of Mathematics, The University of Mississippi,
University, MS 38677, USA. \email{sang@olemiss.edu}}
\author{Lin Ge}\address{Division of Arts and Sciences, Mississippi State University at Meridian,
Meridian, MS 39307, USA. \email{lg481@msstate.edu}}
\date{\today}
\begin{abstract} 
In this paper, we study the self-normalized Cram\a'{e}r-type moderate  deviations for centered independent random variables $X_1, X_2,...$ with $0<E |X_i|^3 <\infty$. The main results refine  Theorems 1.1 and 1.2 of Wang \cite{Wang}, the Berry-Esseen bound (2.11) and Corollaries 2.2 and 2.3 of Jing, Shao and Wang \cite{JingShaoWang} under stronger moment conditions. 
\end{abstract}
\begin{resume} 
\end{resume}
\maketitle
\section*{Introduction}
Let $X_1, X_2, ...$ be independent random variables with $EX_i=0$ and $0<EX_i^2 <\infty$. Set
\begin{gather*}
S_n = \sum_{i=1}^n X_i, \quad V_n^2 = \sum_{i=1}^n X_i^2, \  \   \textrm{and} \ \ B_n^2 = \sum_{i=1}^n EX_i^2.
\end{gather*}

The last two decades have witnessed a significant development on the limit theorems in the self-normalized form $S_n/V_n$,  including central limit theorem, weak invariance principle, law of the iterated logarithm, Berry-Esseen inequality,  large and moderate deviation probabilities. The last two in this list are the main approaches for estimating the error of the normal approximation of the self-normalized probabilities. One advantage of these self-normalized limit theorems is that they usually require less moment conditions than those for the corresponding regular limit theorems. An incomplete list of reference includes Griffin and Kuelbs \cite{Griffin1, Griffin2},  Gin\'{e}, G\"{o}tze and Mason \cite{GineGotzeMason},  
Shao \cite{Shao,Shao3},  Wang and Jing \cite{WangJing}, Cs\"{o}rg\H{o}, Szyszkowicz and Wang \cite{CsorgoSzyszkowiczWang}, Jing, Shao and Wang \cite{JingShaoWang}, Jing, Shao and Zhou \cite{Jing2}, Robinson and Wang \cite{RobinsonWang}, Wang \cite{Wang05, Wang}, the survey papers of Shao \cite{Shao1, Shao2} and Shao and Wang \cite{ShaoWang}.  A systematic treatment of self-normalized limit theory is also collected in the book by de la Pe\~{n}a, Lai and Shao \cite{delaPenaLaiShao}.  
 
The focus of this paper is on the self-normalized Cram\a'{e}r-type moderate deviations. Let $b= x/B_n$ and $\tau = B_n/\max \{ 1, x \}$ with $x \geq 0$. Set 
\begin{gather}
L_{kn}= B_n^{-k}\sum_{i=1}^n E X_i^k \  \ \textrm{and} \  \ \bar{L}_{kn}= B_n^{-k}\sum_{i=1}^n E X_i^k I(|bX_i| \leq 1),\ \ k \geq 2, \nonumber\\
\mathcal{L}_{kn}= B_n^{-k}\sum_{i=1}^n E |X_i|^k  \ \ \textrm{and}  \ \  \mathcal{\bar{L}}_{kn}= B_n^{-k}\sum_{i=1}^n E |X_i|^k I(b|X_i| \leq 1), \ \ k\geq 2,\nonumber\\
\Delta^{j,k}_{n,x}=\tau^{-j}\sum_{i=1}^n E|X_i|^j I (|X_i|  > \tau)+ \tau^{-k}\sum_{i=1}^n E|X_i|^{k}I(|X_i| \leq \tau), \ \   k >j\geq 2. \nonumber
\end{gather}

For the i.i.d. case, Shao \cite{Shao3} proved that as $n \rightarrow \infty$,
\begin{eqnarray*}\label{selfnormalizedM}
\frac{P (S_n/V_n \geq x)}{1-\Phi(x)}\rightarrow 1,   \;\;  \ \  \frac{P (S_n/V_n \leq -x)}{\Phi(-x)}\rightarrow 1
\end{eqnarray*}
holds uniformly for $x\in [0, o(n^{\delta/(4+2\delta)})]$ under the conditions $\mathbb{E}X_1=0$ and $\mathbb{E}|X_1|^{2+\delta}<\infty$ for $0<\delta\le 1$. Here  $\Phi (x)$ is the distribution function of the standard normal random variables.

For independent random variables $X_1, X_2, \dots$ with $EX_i=0$, $0< EX_i^2<\infty$,  $x^2 \max_{i}EX_i^2 \leq B_n^2$, and $ \Delta^{2,3}_{n,x}\leq (1+x)^2/A$ where $A$ is a constant sufficiently large,  Jing, Shao and Wang \cite{JingShaoWang} established a Cram\a'{e}r-type deviation result for self-normalized sums
\begin{eqnarray*}
\frac{ P (S_n/V_n \geq x)}{1-\Phi(x)}=e^{O(1)\Delta^{2,3}_{n,x}}.
\end{eqnarray*}

Wang \cite{Wang} developed new techniques to refine the self-normalized Cram\a'{e}r-type deviation results and simplified the proof. Theorem 1.2 there states that if $EX_i=0$, $0< E |X_i|^3< \infty$, $x \mathcal{L}_{3n} \leq 1/A$, $x^3  \max_{i}E|X_i|^3 \leq B^3_n/27$, and $|c| \leq x/5$, then
\begin{eqnarray}
P(S_n \geq x V_n +c B_n) = \tilde{\Psi}_x (\lambda_1, c)(1-\Phi(x+c)) e^{O(1) \Delta^{3,4}_{n,x}}\{1+O(1)(1+x)\mathcal{L}_{3n}\}, \nonumber
\end{eqnarray}
where $\tilde{\Psi}_x (t, c)= e^{  m(t)+(x+c)^2/2-t(x^2+2xc) }$ for $m(t)=\sum_{i=1}^n \log E e^{t (2 bX_i-b^2 X_i^2)}$, and $\lambda_1$ is the solution of $\lambda$ for $m'(\lambda)=x^2 +2cx$. His Theorem 1.1 states that if $EX_i=0$, $0<EX_i^4<\infty$ and $x \leq \mathcal{L}_{4n}^{-1/4}$, then
\begin{eqnarray}
P(S_n/V_n \geq x  ) = \tilde{\Psi}_x (\lambda_0, 0)(1-\Phi(x)) \{1+O(1)(1+x)\mathcal{L}_{3n}+O(1) (1+x^4)\mathcal{L}_{4n}\},\nonumber
\end{eqnarray}
where $\lambda_0$ is the solution of $\lambda$ for $m'(\lambda)=x^2$. 

In this paper, we further refine the proofs and results of Wang \cite{Wang}. Section \ref{Section2} gives the main results including a theorem and three corollaries. The theorem refines Theorems 1.1 and 1.2 of Wang \cite{Wang} under stronger moment conditions. The corollaries refine the Berry-Esseen bound (2.11) and Corollaries 2.2 and 2.3 of Jing, Shao and Wang \cite{JingShaoWang} under stronger moment conditions. The proof of the theorem  is also given in Section 2 using the propositions in Section \ref{Section3}. In Proposition \ref{corollary1}, we obtain a formula for the probability $P \{2bS_{n} - (\alpha_0-\alpha_3 x \bar{L}_{3n})b^2 V_{n}^2\geq (2-\alpha_0)x^2 +\delta(x) \}$ where $\alpha_0$ and $\alpha_3$ are constants with different values in different situations, and $\delta(x)$ is a function of $x$. This probability generalizes the probability $P(2bS_n -b^2 V_n^2 \geq x^2)$ given by Jing, Shao and Wang \cite{JingShaoWang} and the probability $P\{2bS_n -b^2 V_n^2 \geq x^2 +\delta(x)\}$ given by Wang \cite{Wang}. This generalized probability is necessary to produce the desired accuracy in Theorem \ref{Thm1}, in particular, the exponent $-x^3 L_{3n}/3-x^4 L_{4n}/12$ in Equations (\ref{Eq4}) and (\ref{Eq5}).



\section {Main Results}\label {Section2}

\begin{theorem}\label {Thm1}
Let $X_1, X_2, ...$ be independent random variables with $EX_i=0$ and $0<E|X_i|^3 <\infty$. Assume that there exists a constant $1<A<\infty$ sufficiently large such that for $x \geq 0$,
\begin{eqnarray} \label {condition1}
 x\mathcal{L}_{3n}\leq 1/A  
\end{eqnarray} 
and
\begin{eqnarray} \label {condition2}
 x^2 B_n^{-2} \max_{1 \leq i \leq n} EX_i^2 \leq  1/12.
\end{eqnarray}
Then
\begin{eqnarray}\label{Eq1}
\frac{P(S_n \geq x V_n) }{1-\Phi(x) }&=& \exp \left ( - \frac{x^3 \bar{L}_{3n}}{3 }  - \frac{x^4 \bar{ L}_{4n}}{12 }+O(1) \Delta^{3,5}_{n,x}\right )  \{1+ O(1) (1+x) \mathcal{L}_{3n}   \}.
\end{eqnarray}
Consequently, 
\begin{eqnarray}\label {Eq3}
\frac{P(S_n \geq x V_n)}{1-\Phi(x) } = \exp \left ( - \frac{x^3 L_{3n}}{3 }+ O(1)\Delta^{3,4}_{n,x} \right) \Big \{1+O(1)(1+x) \mathcal{L}_{3n}\Big \}.
\end{eqnarray}
If $0< EX_i^4 < \infty$, then 
\begin{eqnarray}\label{Eq4}
\frac{P(S_n \geq x V_n) }{1-\Phi(x)}&=&  \exp \left ( - \frac{x^3 L_{3n}}{3 }-\frac{x^4 L_{4n}}{12 }+O(1) \Delta^{4,5}_{n,x} \right ) \{1+O(1)  (1+x) \mathcal{L}_{3n}   \}.
\end{eqnarray}
If $0< E|X_i|^{4+\delta} < \infty$ with $0<\delta\leq 1$, then for $0 \leq x \leq \mathcal{L}_{(4+\delta)n}^{-1/(4+\delta)}$,
\begin{eqnarray}\label{Eq5}
\frac{P(S_n \geq x V_n) }{1-\Phi(x)}&=&  \exp \left ( - \frac{x^3 L_{3n}}{3 }-\frac{x^4 L_{4n}}{12 } \right ) \big \{1+O (1) (1+x) \mathcal{L}_{3n}+O(1)( 1+x)^{4+\delta} \mathcal{L}_{(4+\delta)n}\big \}.
\end{eqnarray}
\end{theorem}
\begin{proof}
First we prove the lower bound. By Lemma \ref{properties}(i), we have $x^3 \mathcal{\bar{L}}^3_{3n}  \leq  x^3\mathcal{\bar{L}}_{5n} \leq \Delta_{n,x}^{3,5}/x^2$. Since $|(1-s)^{-1}-1-s-s^2| \leq 2 |s|^3$ for $|s| \leq 1/2$, then 
\begin{eqnarray}\label {DC109}
 \bigg|  \frac{1}{1- x \bar{L}_{3n}/2 } -1 -\frac{ x \bar{L}_{3n}}{2} -\frac{ x^2  \bar{L}^2_{3n}}{4}\bigg| \leq \frac{x^3  \mathcal{\bar{L}}^3_{3n} }{4} \leq \frac{\Delta^{3,5}_{n,x}}{4x^2}. 
\end{eqnarray}
Hence
\begin{eqnarray}\label {BN25}
&& P(S_n \geq x V_n)\nonumber\\
 &=& P \bigg \{  2bS_n - (1- x \bar{L}_{3n}/2) b^2 V_n^2 \geq \frac{x^2}{1- x \bar{L}_{3n}/2}\nonumber \\
&&- \bigg ( \frac{x}{  (1- x \bar{L}_{3n}/2)^{1/2}} -b V_n  (1- x \bar{L}_{3n}/2)^{1/2}  \bigg )^2 \bigg \} \nonumber\\
&\geq& P \bigg \{ 2bS_n - (1 - x \bar{L}_{3n}/2) b^2 V_n^2 \geq \frac{x^2}{1- x \bar{L}_{3n}/2}  \bigg \}\nonumber\\
&\geq& P \bigg \{ 2bS_n - (1- x \bar{L}_{3n}/2) b^2 V_n^2 \geq x^2 +\frac{x^3 \bar{L}_{3n}}{2} + \frac{x^4 \bar{L}^2_{3n}}{4} +\Delta^{3,5}_{n,x}\bigg \} \nonumber.
\end{eqnarray}
By Proposition \ref{corollary1} with $\alpha_0=1, \ \alpha_3 =1/2, \  \beta_3 = 1/2, \ \beta_4 = 0$ and $\beta_5=1/4$,
\begin{eqnarray}\label {8L}
&& P(S_n \geq x V_n) \nonumber\\
&\geq& \exp \bigg \{  -\frac{1}{3 }x^3 \bar{L}_{3n} -\frac{1}{12}x^4 \bar{L}_{4n}-A \Delta^{3,5}_{n,x}\bigg \} (1-\Phi(x))(1-Ax \mathcal{L}_{3n}).  \nonumber
\end{eqnarray}

Next we prove the upper bound. First we consider $0 \leq x \leq 2$. Observe that $ x^4 \bar{L}_{4n} \leq x^3 \mathcal{\bar{L}}_{3n} \leq x^2 /A \leq 1$ by condition (\ref{condition1}). Since $|e^s - 1| \leq e^{s \vee 0}|s| $, then 
\begin{eqnarray}
\Big |\exp  \Big (-\frac{1}{3} x^3 \bar{L}_{3n}-\frac{1}{12}x^4\bar{L}_{4n}\Big) -1 \Big | \leq \Big (\frac{1}{3} +\frac{1}{12}\Big )e^{1/3} x^3 \mathcal{\bar{L}}_{3n} \leq 3 x \mathcal{\bar{L}}_{3n}. \nonumber
\end{eqnarray}
Hence by condition (\ref{condition1}),
\begin{eqnarray} \label {5A}
&&\exp  \Big (-\frac{1}{3} x^3 \bar{L}_{3n}-\frac{1}{12}x^4\bar{L}_{4n}\Big) \{1+ 2A (1+x) \mathcal{\bar{L}}_{3n}\} \nonumber\\
&\geq& (1-3 x \mathcal{\bar{L}}_{3n}) \{1+ 2A (1+x) \mathcal{\bar{L}}_{3n}\} \geq   1+   A(1+x) \mathcal{\bar{L}}_{3n}  .
\end{eqnarray}
By (2.16) of Wang \cite{Wang}, $|P(S_n \geq x V_n) -( 1-\Phi(x))| \leq A  \mathcal{L}_{3n}$. Then for $0 \leq x \leq 2$,
\begin{eqnarray}\label {6A}
\frac{P(S_n \geq x V_n)}{1-\Phi(x)} \leq 1 +A(1+x) \mathcal{L}_{3n}.
\end{eqnarray}
Combining (\ref{5A}) and (\ref{6A}), we obtain the upper bound. For $x \geq 2$,
\begin{eqnarray}
P(S_n \geq x V_n)= P\Big (S_n \geq x V_n, \  \max_{1\leq i \leq n}|X_i|>\tau \Big) + P(\bar{S}_n \geq x \bar{V}_n) \nonumber
\end{eqnarray}
where $\tau=B_n/x$, $\bar{S}_n=\sum_{i=1}^n X_i I(|X_i| \leq \tau)$ and $\bar{V}^2_n = \sum_{i=1}^n X^2_i I(|X_i| \leq \tau)$. 
Page 2181 of Jing, Shao and Wang \cite{JingShaoWang} shows that
\begin{eqnarray}\label {DD120}
P\Big (S_n \geq x V_n, \  \max_{1\leq i \leq n}|X_i|>\tau \Big) &\leq& \sum_{i=1}^n P \left (S_n^{(i)} \geq (x^2-1)^{1/2} V^{(i)}_n \right) P(|X_i| > \tau).\nonumber
\end{eqnarray}
Since $\sum_{i=1}^nP(|X_i| > \tau) \leq \Delta^{3,5}_{n,x}$, then the upper bound follows from Propositions \ref{Proposition2} and \ref{Prop1}. 
\end{proof}
\bigskip

Note that Equations (\ref{Eq4}) and (\ref{Eq5}) refine Theorems 1.2 and 1.1 (when $c=0$) of Wang \cite{Wang}, respectively, under higher moment conditions.

\begin{corollary}\label {Corollary1}
Suppose that $\max_{1 \leq i \leq n}EX_i^4 \leq C$ for some $C <\infty$ and $ B_n^2  \geq cn$ for some $c>0$. If $\sum_{i=1}^n \mathbb{E}X_i^3 = O(n^{\gamma})$ for $0 \leq \gamma \leq 1$, then 
\begin{eqnarray}
P (S_n \geq xV_n) - (1-\Phi(x) )=O\left (\frac{x^2}{n^{3/2-\gamma}} + \frac{x^3}{n}+\frac{1}{\sqrt{n}} \right ) e^{-x^2/2} \nonumber
\end{eqnarray}
for $0 \leq x \leq \min \{O(n^{1/2-\gamma/3}), \  O(n^{1/4}) \}$.
\end{corollary}
\begin{proof}
Let $y=-x^3 L_{3n}/3 -x^4 L_{4n}/12+O(1)\Delta^{4,5}_{n,x}$ in (\ref{Eq4}). Note that $\Delta^{4,5}_{n,x} \leq (1+x)^4 L_{4n}$. Then $y$ is bounded for $x = \min \{O(n^{1/2-\gamma/3}), \  O(n^{1/4}) \}$. Since $|e^y-1| \leq e^{y \vee 0}|y|$, then $e^y= 1+ O(1)y$. Hence $P(S_n \geq x V_n)=(1-\Phi(x))\{1+O(1)y+O(1)(1+x )\mathcal{L}_{3n}\}$ by (\ref{Eq4}). Note that $1-\Phi(x)\leq 2 e^{-x^2/2}/(1+x)$ for $x\geq0$. Then $x^3 L_{3n}(1-\Phi(x))= O(x^2/n^{3/2-\gamma})e^{-x^2/2}$ and $x^4 L_{4n}(1-\Phi(x))= O(x^3/n )e^{-x^2/2}$. Moreover, $(1+x)\mathcal{L}_{3n}(1-\Phi(x))=O(1/\sqrt{n})e^{-x^2/2}$ because $E|X_i|^3 \leq (EX_i^4)^{3/4} \leq C^{3/4}$. Therefore, the corollary follows. 
\end{proof}
\bigskip

Compared with the Berry-Esseen bound (2.11) by Jing, Shao and Wang \cite{JingShaoWang}, this corollary shows that the bound can be lowered and the range of $x$ can be extended under stronger moment conditions. For example, if $\max_{1 \leq i \leq n} E|X_i|^3 \leq C$ and $B^2_n \geq cn$, their result shows that $P (S_n > xV_n) - (1-\Phi(x)) = O(1)(1+x)^2e^{-x^2/2}/\sqrt{n}$ for $x =O(n^{1/6})$. However, if $\max_{1 \leq i \leq n} EX_i^4\leq C$ and $B_n^2 \geq cn$, the above corollay shows that $P (S_n > xV_n) - (1-\Phi(x))=O(1)(1+x)e^{-x^2/2}/\sqrt{n}$ for $\sum_{i=1}^n \mathbb{E}X_i^3 = O(n^{3/4})$ and $x = O(n^{1/4})$.

\begin{corollary}\label {Corollary2}
Suppose that $\max_{1 \leq i \leq n}EX_i^4 \leq C$ for some $C <\infty$ and $  B_n^2  \geq c n$ for some $c>0$. If $\sum_{i=1}^n \mathbb{E}X_i^3 = O(n^{\gamma})$ for $0 \leq \gamma \leq 1$, then by (\ref{Eq4}),
\begin{eqnarray}\label {sM}
\frac{P (S_n \geq  xV_n)}{1-\Phi(x)} \rightarrow 1 \nonumber
\end{eqnarray}
 uniformly for $0 \leq x \leq \min \{ o(n^{1/2-\gamma/3}), \ o(n^{1/4})\}$.
\end{corollary}

This corollary extends the range of $x$ for Corollary 2.2 of Jing, Shao and Wang \cite{JingShaoWang} where $0 \leq x \leq O(n^{\delta/(4+2 \delta)})$ for $0< \delta<1$.

\begin{corollary} \label {Corollary3}
Let $X_1, X_2, ...$ be i.i.d. random variables with $EX_1=0$, $\sigma^2=E X_1^2$ and $0<E|X_1|^3 <\infty$. Assume that $ 0\leq x  \leq \sqrt{n} /A$ for some sufficiently large constant $2\sqrt{3} \leq A<\infty$.  Then
\begin{eqnarray*}
\frac{P(S_n \geq x V_n)}{1-\Phi(x) } = \exp \left ( - \frac{x^3 E X_1^3}{3 \sqrt{n}\sigma^3 }+ O(1)\Delta^{3,4}_{n,x} \right) \left  \{1+O(1)\frac{1+x}{\sqrt{n}}  \right \},
\end{eqnarray*}
where $\Delta^{3,4}_{n,x}=\frac{(1 \vee x)^3 E |X_1|^3 I\{|X_1|> \sqrt{n}\sigma/(1 \vee x)\}}{\sqrt{n}\sigma^3}+\frac{(1 \vee x)^4 E X_1^4 I\{|X_1|\leq \sqrt{n}\sigma/(1 \vee x)\}}{n \sigma^4} $. \\If $0< EX_1^4 < \infty$, then 
\begin{eqnarray*}
\frac{P(S_n \geq x V_n) }{1-\Phi(x)}&=&  \exp \left ( - \frac{x^3 E X_1^3}{3 \sqrt{n}\sigma^3 }-\frac{x^4 E X_1^4}{12n\sigma^4 }+O(1) \Delta^{4,5}_{n,x} \right )  \left \{1+O(1)  \frac{1+x}{\sqrt{n}} \right \},
\end{eqnarray*}
where $\Delta^{4,5}_{n,x}=\frac{(1 \vee x)^4 E X_1^4 I\{|X_1|> \sqrt{n}\sigma/(1 \vee x)\}}{n\sigma^4}+\frac{(1 \vee x)^5 E |X_1|^5 I\{|X_1|\leq \sqrt{n}\sigma/(1 \vee x)\}}{n^{3/2} \sigma^5} $.\\
If $0< E|X_1|^{4+\delta} < \infty$ with $0<\delta\leq 1$, then for $0 \leq x \leq n^{(2+\delta)/(8+2\delta)}\sigma/(E|X_1|^{4+\delta})^{1/(4+\delta)}$,
\begin{eqnarray*}
\frac{P(S_n \geq x V_n) }{1-\Phi(x)}&=&  \exp \left ( - \frac{x^3 E X_1^3}{3 \sqrt{n}\sigma^3 }-\frac{x^4 E X_1^4}{12n\sigma^4 } \right )  \left \{1+\frac{O(1) (1+x)}{\sqrt{n}}+\frac{O(1)( 1+x)^{4+\delta}}{ n^{1+\delta/2}}\right\}.
\end{eqnarray*}
\end{corollary}
Note that since $X_1, X_2, ...$ are i.i.d. random variables, condition (\ref{condition1}) becomes $ x \leq \sqrt{n}\sigma^3/(AE|X_1|^3)\leq \sqrt{n}/A$, and (\ref{condition2}) becomes $x \leq    \sqrt{n/12}$.


\section{Lemmas and Propositions}\label {Section3}

From now on, all $|O(1)| \leq A$. First we establish some preliminary facts in the following lemma.

\begin{lemma} \label {properties}
(i)
\begin{eqnarray*}\label {porp1}
\mathcal{\bar{L}}^2_{4n}\leq    \mathcal{\bar{L}}_{3n} \mathcal{\bar{L}}_{5n}, \ \  \mathcal{\bar{L}}_{3n} \mathcal{\bar{L}}_{4n} \leq \mathcal{\bar{L}}_{5n} \ \ \textrm{and} \ \  \mathcal{\bar{L}}^3_{3n}  \leq  \mathcal{\bar{L}}_{5n}.
\end{eqnarray*}
(ii)
\begin{eqnarray*}\label {prop2}
x^2 \mathcal{\bar{L}}_{4n} \leq  \max \{ \mathcal{\bar{L}}_{3n},\ x^4 \mathcal{\bar{L}}_{5n} \}.
\end{eqnarray*}
\end{lemma}
\begin{proof}
(i) Let $\bar{X}_i=X_i I(|bX_i| \leq 1)$. Then $E\bar{X}_i^4 \leq (E|\bar{X}_i|^3)^{1/2}(E|\bar{X}_i|^5)^{1/2}$. Hence 
$$ \left ( \sum_{i=1}^n E\bar{X}_i^4 \right )^2 \leq  \sum_{i=1}^n E |\bar{X}_i|^3 \sum_{i=1}^n E|\bar{X}_i|^5$$ 
by the Cauchy-Schwarz inequality. Therefore,
\begin{eqnarray}\label {BP28}
\mathcal{\bar{L}}^2_{4n}\leq    \mathcal{\bar{L}}_{3n} \mathcal{\bar{L}}_{5n}.
\end{eqnarray}
Similarly,
\begin{eqnarray}\label {BP9}
\mathcal{\bar{L}}_{3n}^2 \leq \mathcal{\bar{L}}_{2n}\mathcal{\bar{L}}_{4n} \leq   \mathcal{\bar{L}}_{4n}.
\end{eqnarray}
By (\ref{BP28}) and (\ref{BP9}), $\mathcal{\bar{L}}_{3n} \mathcal{\bar{L}}_{4n} \leq \mathcal{\bar{L}}_{5n}$ which, together with (\ref{BP9}), implies $\mathcal{\bar{L}}^3_{3n} \leq \mathcal{\bar{L}}_{5n}$. 
 
(ii) If $x^2 \mathcal{\bar{L}}_{4n} \leq    x^4 \mathcal{\bar{L}}_{5n}$, then $x^2 \mathcal{\bar{L}}_{4n} \leq \max \{ \mathcal{\bar{L}}_{3n},\ x^4 \mathcal{\bar{L}}_{5n} \}$. If $x^2 \mathcal{\bar{L}}_{4n} \geq    x^4 \mathcal{\bar{L}}_{5n}$, then $ \mathcal{\bar{L}}^2_{4n}\leq     \mathcal{\bar{L}}_{3n} x^{-2}\mathcal{\bar{L}}_{4n} $ by (\ref{BP28}) and hence $x^2 \mathcal{\bar{L}}_{4n} \leq  \mathcal{\bar{L}}_{3n}\leq \max \{ \mathcal{\bar{L}}_{3n},\ x^4 \mathcal{\bar{L}}_{5n} \}$. 
\end{proof}
Let $b= x/B_n$ with $x>0$. Define the function
\begin{eqnarray}\label {definitionxi}
\xi(y)=2b y -\big(\alpha_0-\alpha_3 x \bar{L}_{3n}\big)b^2 y^2 \nonumber
\end{eqnarray}
where $\alpha_0$ and $\alpha_3$ are constants such that $1/3 \leq   \alpha_0-\alpha_3 x \bar{L}_{3n} \leq 7/3$. Here $1/3$ and $7/3$ are selected from the estimates of $I_2$ and $I_3$ in Proposition \ref{Proposition2}. In other situations, we only need $\alpha_0=1$ and $\alpha_3 = 0$ or $1/2$, and thus $\alpha_0-\alpha_3 x \bar{L}_{3n} $ is close to $1$ by condition (\ref{condition1}). 

Denote that $X_{i, (1)} =X_{i}$ and $X_{i, (2)} =X_i I(|bX_i|\leq 1)$ for $i \geq 1$. For $\lambda>0$, let $Z_{\lambda, 1, (k)},..., Z_{\lambda, n, (k)}$ be independent random variables with $Z_{\lambda, i, (k)}$ having the distribution function
\begin{eqnarray}\label {CD30}
P(Z_{\lambda, i, (k)} \leq u) = \frac{\int_{-\infty}^u e^{\lambda \xi(y)} dP(X_{i, (k)} \leq y)}{E e^{\lambda \xi(X_{i, (k)})}} 
\end{eqnarray}
where $k=1$ or $2$. Then for any function $f$,
\begin{eqnarray}\label {Efz}
E f(Z_{\lambda, i, (k)}) = \frac{\int_{-\infty}^{\infty} f(u) e^{\lambda \xi(u)} d P (X_{i, (k)} \leq u)}{E e^{\lambda \xi(X_{i, (k)})}}.
\end{eqnarray}
In particular,
\begin{eqnarray}\label {CD32}
E  \xi (Z_{\lambda, i, (k)}) = \frac{\int_{-\infty}^{\infty}\xi(u) e^{\lambda \xi(u)} d P (X_{i, (k)} \leq u)}{E e^{\lambda \xi(X_{i, (k)})}}  = \frac{d \big ( \log E e^{\lambda \xi(X_{i, (k)})} \big )}{d \lambda}
\end{eqnarray}
and 
\begin{eqnarray}\label {CD32a}
Var (\xi(Z_{\lambda,i, (k)}))= \frac{d^2 \big ( \log E e^{\lambda \xi(X_{i, (k)})} \big)}{d \lambda^2}.
\end{eqnarray}

\begin{lemma} \label {Lemma2}
Let $b=x/B_n$ with $x>0$, $X_{i, (1)}=X_i$ and $X_{i, (2)}= X_i I(|bX_i| \leq 1)$. Under conditions (\ref{condition1}) and (\ref{condition2}), for $7/16 \leq \lambda \leq 9/16$ and $k=1$ or $2$,
\begin{eqnarray} \label {noxi}
\sum_{i=1}^n \log E e^{\lambda \xi(X_{i, (k)})}&=&(2 \lambda^2-\lambda \alpha_0)x^2  + \left (\frac{4\lambda^3}{3} -2 \lambda^2 \alpha_0 +\lambda \alpha_3 \right )x^3 \bar{L}_{3n} \nonumber\\
&& + \left ( \frac{2\lambda^4}{3}-2 \lambda^3 \alpha_0+\frac{\lambda^2 \alpha_0^2}{2} \right )x^4 \bar{L}_{4n} +2 \lambda^2 \alpha_3 x^4 \bar{L}^2_{3n}  \nonumber\\
&&-\frac{\left(2 \lambda^2-\lambda \alpha_0 \right)^2}{2} \frac{x^4\sum_{i=1}^n ( E \bar{X}_i^2)^2}{B_n^4} +O(1)\Delta^{3,5}_{n,x}.
\end{eqnarray}
Consequently,
\begin{eqnarray} \label {Eeta1}
\sum_{i=1}^n E \xi (Z_{\lambda, i, (k)}) &=&(4 \lambda -\alpha_0)  x^2  +(4 \lambda^2-4\lambda \alpha_0 +\alpha_3) x^3 \bar{L}_{3n}  \nonumber\\
&&+\left ( \frac{8\lambda^3}{3}-6\lambda^2 \alpha_0+\lambda \alpha_0^2\right ) x^4 \bar{L}_{4n} +4 \lambda \alpha_3 x^4 \bar{L}^2_{3n} \nonumber\\
&& +\left(-8\lambda^3+6 \lambda^2 \alpha_0- \lambda \alpha_0^2 \right) \frac{x^4\sum_{i=1}^n ( E \bar{X}_i^2)^2}{B_n^4} \nonumber\\
&&+O(1) \Delta^{3,5}_{n,x} \qquad
\end{eqnarray}
and  
\begin{eqnarray} \label {Vareta1}
\sum_{i=1}^n Var( \xi (Z_{\lambda, i, (k)}))= 4 x^2 +O(1)  x^3 \mathcal{L}_{3n}.
\end{eqnarray}
Moreover,  for $  m \geq 2$ and $|O(1)|\leq A^{1/5}$,
\begin{eqnarray} \label {CD31}
E Z_{\lambda, i, (2)}^m= E  \bar{X}_{i}^m + 2 \lambda b E \bar{X}_i^{m+1} +O(1) b^2 E| \bar{X}_{i}|^{m+2}.
\end{eqnarray}
\end{lemma}
\begin{proof} 
We follow some proofs in Lemmas 6.1 and 6.2 of Jing, Shao and Wang \cite{JingShaoWang}. Let $\gamma = 2 \lambda$ and $\theta = \lambda  (\alpha_0 - \alpha_3 x \bar{L}_{3n})$. Then $\lambda \xi(X_i)= \gamma b X_i - \theta b^2 X_i^2$. Hence
\begin{eqnarray}\label {A29a}
E (e^{\lambda \xi(X_{i, (k)})}-1) I(|bX_i|>1) \leq   e^{\gamma^2/(4 \theta)} b E|X_i|  I(|bX_i|>1).
\end{eqnarray}
Since $|e^s-1-s-s^2/2-s^3/6-s^4/24| \leq |s|^5 e^{s \vee 0}/120$ for $s \in R$, then
\begin{eqnarray}\label {A32a}
&& E (e^{\lambda \xi (X_{i, (k)})}-1) I(|bX_i|\leq 1)\nonumber\\
&=&   E \left(\gamma b \bar{X}_i - \theta b^2 \bar{X}_i^2 \right) + \frac{1}{2}E \left(\gamma b \bar{X}_i - \theta b^2 \bar{X}_i^2 \right)^2+ \frac{1}{6} E \left(\gamma b \bar{X}_i - \theta b^2 \bar{X}_i^2 \right)^3\nonumber\\
&& + \frac{1}{24} E \left(\gamma b \bar{X}_i - \theta b^2 \bar{X}_i^2 \right)^4+O(1) \frac{e^{\gamma^2/(4\theta)}}{120}E \left |\gamma b\bar{X}_i - \theta b^2 \bar{X}_i^2 \right|^5  \nonumber\\
&=& \gamma b E \bar{X}_i + \left(\frac{\gamma^2}{2}-\theta \right)b^2 E \bar{X}_i^2 + \left (\frac{\gamma^3}{6} -\gamma \theta\right )b^3 E \bar{X}_i^3 \nonumber \\
&&+ \left (\frac{\gamma^4}{24}-\frac{\gamma^2 \theta}{2}+  \frac{\theta^2}{2}\right )b^4 E \bar{X}_i^4 + O'_{\gamma, \theta}  b^5 E|\bar{X}_i|^5 ,
\end{eqnarray}
where 
\begin{eqnarray}\label {21kg}
|O'_{\gamma, \theta}| \leq \left \{ \frac{1}{6}(\gamma+\theta)^3 +\frac{1}{24}(\gamma + \theta)^4 +\frac{e^{\gamma^2/(4 \theta)}}{120}(\gamma+\theta)^5  \right\}.
\end{eqnarray}
By (\ref{A29a}) and (\ref{A32a}),
\begin{eqnarray}\label {CC20}
E e^{\lambda \xi(X_{i,(k)})}&=&1 + \gamma b E\bar{X}_i+\left(\frac{\gamma^2}{2}-\theta \right)b^2 E \bar{X}_i^2 + \left (\frac{\gamma^3}{6} -\gamma \theta   \right )b^3 E \bar{X}_i^3 \nonumber\\
&&+ \left ( \frac{\gamma^4}{24}-\frac{\gamma^2 \theta}{2}+ \frac{\theta^2}{2}\right )b^4 E \bar{X}_i^4 \nonumber\\
&& +O(1)  e^{\gamma^2/(4 \theta)} b E|X_i| I(|bX_i|>1)+O_{\gamma, \theta}' b^5 E |\bar{X}_i|^5 ,\qquad
\end{eqnarray}
where $|O(1)| \leq 1$. We want to show that $|E e^{\lambda \xi(X_{i, (k)})}-1| \leq 11/12$. Since $7/16 \leq \lambda \leq 9/16$ and $1/3 \leq   \alpha_0-\alpha_3 x \bar{L}_{3n}\leq 7/3$, then $7/8 \leq \gamma \leq 9/8$ and $7/48 \leq \theta \leq 63/48$. Hence by (\ref{21kg}),
\begin{eqnarray}
&& \left |\gamma^2/2-\theta \right |+ \left |\frac{\gamma^3}{6} -\gamma \theta   \right | + \left | \frac{\gamma^4}{24}-\frac{\gamma^2 \theta}{2}+ \frac{\theta^2}{2}\right | + |O'_{\gamma, \theta}|\nonumber\\
& \leq&   \left | (7/8)^2/2- 63/48|+|  (7/8)^3/6 -(9/8)(63/48)\right |\nonumber\\
&&+\left |(9/8)^4/24-(7/8)^2(7/48)/2+(63/48)^2/2 \right| \nonumber\\
&&+ \left |(9/8+63/48)^3/6 +(9/8+63/48)^4/24+e^{27/16}(9/8+63/48)^5/120 \right| \nonumber\\
&\leq&0.93+ 1.37 + 0.88+7.77 = 10.95.\nonumber
\end{eqnarray}
Note that $\gamma b E \bar{X}_i = - \gamma b E X_i I(|bX_i|>1)$. Since $\gamma^2/(4\theta) = \lambda/(\alpha_0 - \alpha_3 x \bar{L}_{3n}) \leq 27/16$, then
\begin{eqnarray}
\big |\gamma +e^{\gamma^2/(4 \theta)}\big|   b E|X_i| I(|bX_i|>1) &\leq&  (9/8+e^{27/16}) b^2 E X_i^2 I(|bX_i|>1) \nonumber\\
&\leq& 6.54b^2 E X_i^2 I(|bX_i|>1) .\nonumber
\end{eqnarray}
From the above, $|E e^{\lambda \xi(X_{i, (k)})}-1| \leq 10.95 b^2 EX_i^2 \leq 11/12$ by condition (\ref{condition2}). Note that $|\log (1+y)-y+y^2/2| \leq 12^3 |y|^3/3$ for $|y| \leq 11/12$. Since $\gamma = 2 \lambda$ and $\theta = \lambda  (\alpha_0- \alpha_3 x \bar{L}_{3n})$, then by (\ref{CC20}) and Lemma \ref{properties}(i),
\begin{eqnarray}\label {32ab}
&&\log E e^{\lambda \xi(X_{i, (k)})} =\log \left (1+ \{ E e^{\lambda \xi(X_{i, (k)})} -1 \}\right)\nonumber\\
&=&\left \{2 \lambda^2-\lambda  (\alpha_0- \alpha_3 x \bar{L}_{3n})\right \}b^2 E \bar{X}_i^2 + \left \{\frac{4\lambda^3}{3} -2 \lambda^2  (\alpha_0- \alpha_3 x \bar{L}_{3n}) \right \}b^3 E \bar{X}_i^3 \nonumber \\
&&+ \left \{ \frac{2\lambda^4}{3}-2 \lambda^3  (\alpha_0- \alpha_3 x \bar{L}_{3n})+\frac{\lambda^2  (\alpha_0- \alpha_3 x \bar{L}_{3n})^2}{2} \right \} b^4  E \bar{X}_i^4 \nonumber\\
&& - \frac{1}{2}\left \{2 \lambda^2-\lambda  (\alpha_0- \alpha_3 x \bar{L}_{3n})\right \}^2 b^4 ( E \bar{X}_i^2)^2  + O(1)b^3 E| X_i|^3 I(|bX_i|>1) \nonumber\\
&&+O(1) b^5E |\bar{X}_i|^5. \nonumber
\end{eqnarray}
Hence (\ref{noxi}) follows from Lemma \ref{properties}(i). Taking derivative with respect to $\lambda$ on both sides of (\ref{noxi}) and by (\ref{CD32}), we obtain (\ref{Eeta1}). Taking derivative with respect to $\lambda$ on both sides of (\ref{Eeta1}) and by (\ref{CD32a}), we obtain (\ref{Vareta1}). Similar to (\ref{CC20}),
\begin{eqnarray} \label {CC26}
E \bar{X}_i^m e^{\lambda \xi(X_{i,(2)})}=E\bar{X}_i^m+\gamma b E \bar{X}_i^{m+1} +O(1)  b^2 E |\bar{X}_i|^{m+2}
\end{eqnarray}
for $m \geq 2$. Since $|(1+y)^{-1}-1| \leq 12^2|y|$ for $|y| \leq 11/12$ and $\gamma=2\lambda$, then by (\ref{Efz}), (\ref{CC20}) and (\ref{CC26}), 
\begin{eqnarray*}
E  Z_{\lambda, i, (2)}^m= \frac{E \bar{X}_i^m e^{\lambda \xi(\bar{X}_{i,(2)})}}{E e^{\lambda \xi(X_{i, (2)})}}= E \bar{X}_i^m + 2 \lambda b E \bar{X}_i^{m+1} +O(1)  b^2 E|\bar{X}_i|^{m+2}.
\end{eqnarray*}
\end{proof}

\begin{lemma}\label {Prop1}
Let $x>0$, $|\delta(x)| \leq x^2/5$ and $k=1$ or $2$. Under conditions (\ref{condition1}) and (\ref{condition2}), for $\sum_{i=1}^n E  \xi (Z_{\lambda, i, (k)})$ in (\ref{Eeta1}), the equation
\begin{eqnarray}\label {CG39}
\sum_{i=1}^n E  \xi (Z_{\lambda, i, (k)})=(2-\alpha_0) x^2 +\delta(x)
\end{eqnarray}
has a unique solution of $\lambda$, denoted as $\lambda_{\delta}$, with $7/16 \leq \lambda_{\delta} \leq 9/16$. If $\delta(x)=\beta_3 x^3 \bar{L}_{3n}+\beta_4 x^4 \bar{L}_{4n}+ \beta_5 x^4 \bar{L}^2_{3n}+O(1)\Delta^{3,5}_{n,x}$ with $|\beta_j|\leq \sqrt{A}$, then 
\begin{eqnarray}\label {CL42}
\lambda_{\delta} &= &\frac{1}{2} + \left (\frac{\beta_3}{4}+ \frac{\alpha_0}{2} - \frac{\alpha_3}{4} -\frac{1}{4}\right )     x \bar{L}_{3n} +  \left ( \frac{ \beta_4}{4} + \frac{3\alpha_0}{8}-\frac{\alpha_0^2}{8}-\frac{1}{12}\right ) x^2 \bar{L}_{4n} \nonumber\\
&&+\left \{\frac{ \beta_5 }{4} -\frac{\alpha_3}{2}+ (\alpha_0-1) \left ( \frac{\beta_3}{4}+\frac{ \alpha_0}{2}-\frac{\alpha_3}{4}-\frac{1}{4} \right)\right \} x^2 \bar{L}^2_{3n} \nonumber\\
&& +\left (\frac{1}{4}-\frac{3 \alpha_0}{8}+\frac{\alpha_0^2}{8}  \right ) \frac{x^2\sum_{i=1}^n ( E \bar{X}_i^2)^2}{B_n^4}  + \frac{O(1)\Delta^{3,5}_{n,x}}{x^2} .
\end{eqnarray}
\end{lemma}
\begin{proof}
Since the left-hand side of (\ref{Vareta1}) is equal to $d(\sum_{i=1}^n E \xi (Z_{\lambda, i, (k)}))/d \lambda$ and the right-hand side is positive, then $\sum_{i=1}^n E \xi (Z_{\lambda, i})$ increases strictly as $\lambda$ increases. By (\ref{Eeta1}), condition (\ref{condition1}) and $|\delta(x)| \leq x^2/5$,
\begin{eqnarray*}
\sum_{i=1}^n E\xi (Z_{7/16, i, (k)}) \leq (2-\alpha_0) x^2 +\delta(x) \leq \sum_{i=1}^n E  \xi (Z_{9/16, i, (k)}).
\end{eqnarray*}
Hence (\ref{CG39}) has a unique soluion $\lambda_{\delta}$ and $7/16 \leq \lambda_{\delta} \leq 9/16$. By (\ref{Eeta1}) and (\ref{CG39}), 
\begin{eqnarray} \label {CG40}
\lambda_{\delta} &=& \frac{1}{2}  + \left (-  \lambda_{\delta}^2+\lambda_{\delta} \alpha_0 -\frac{\alpha_3}{4} \right) x \bar{L}_{3n}   +\left ( -\frac{2\lambda_{\delta}^3}{3}+ \frac{3\lambda_{\delta}^2 \alpha_0}{2}-\frac{\lambda_{\delta} \alpha_0^2}{4}\right ) x^2 \bar{L}_{4n} \nonumber\\
&& - \lambda_{\delta} \alpha_3 x^2 \bar{L}^2_{3n} +\left(2\lambda_{\delta}^3-\frac{3 \lambda_{\delta}^2 \alpha_0}{2}+ \frac{\lambda_{\delta} \alpha_0^2 }{4}\right) \frac{x^2\sum_{i=1}^n ( E \bar{X}_i^2)^2}{B_n^4}  + \frac{O(1)\Delta^{3,5}_{n,x}}{x^2}\nonumber\\
&&+ \frac{\delta(x)}{4x^2}. \nonumber
\end{eqnarray}
If $\delta(x)=\beta_3 x^3 \bar{L}_{3n}+\beta_4 x^4 \bar{L}_{4n}+ \beta_5 x^4 \bar{L}^2_{3n}+O(1)\Delta^{3,5}_{n,x}$, then
\begin{eqnarray} \label {CL46}
\lambda_{\delta} &=& \frac{1}{2}  + \left (\frac{\beta_3}{4}-  \lambda_{\delta}^2+\lambda_{\delta} \alpha_0 -\frac{\alpha_3}{4} \right) x \bar{L}_{3n}  \nonumber\\\
&& +\left (\frac{\beta_4}{4} -\frac{2\lambda_{\delta}^3}{3}+ \frac{3\lambda_{\delta}^2 \alpha_0}{2}-\frac{\lambda_{\delta} \alpha_0^2}{4}\right ) x^2 \bar{L}_{4n}  + \left (\frac{\beta_5}{4}-\lambda_{\delta} \alpha_3 \right )x^2 \bar{L}^2_{3n} \nonumber\\
&&+\left(2\lambda_{\delta}^3-\frac{3 \lambda_{\delta}^2 \alpha_0}{2}+ \frac{\lambda_{\delta} \alpha_0^2 }{4}\right) \frac{x^2\sum_{i=1}^n ( E \bar{X}_i^2)^2}{B_n^4}  +\frac{O(1) \Delta^{3,5}_{n,x}}{x^2}. 
\end{eqnarray}
By Lemma \ref{properties}(i) and (\ref{CL46}), 
\begin{eqnarray*}\label {40aL}
\lambda^2_{\delta} x \bar{L}_{3n}= \frac{1}{4} x \bar{L}_{3n} + \left (\frac{\beta_3}{4}-  \frac{1}{4}+\frac{ \alpha_0}{2} -\frac{\alpha_3}{4} \right) x^2 \bar{L}^2_{3n}  + \frac{O(1)\Delta^{3,5}_{n,x}}{x^2}
\end{eqnarray*}
and
\begin{eqnarray}\label {43aL}
\lambda_{\delta}\alpha_0 x \bar{L}_{3n}&=& \frac{\alpha_0}{2} x \bar{L}_{3n} + \alpha_0\left (\frac{\beta_3}{4}-  \frac{1}{4}+\frac{ \alpha_0}{2} -\frac{\alpha_3}{4} \right) x^2 \bar{L}^2_{3n}  +\frac{O(1)\Delta^{3,5}_{n,x}}{x^2}.\nonumber
\end{eqnarray}
Hence the second term on the right-hand side of (\ref{CL46})
\begin{eqnarray*}
&&\left (\frac{\beta_3}{4}-  \lambda_{\delta}^2+\lambda_{\delta} \alpha_0 -\frac{\alpha_3}{4} \right) x \bar{L}_{3n} \nonumber\\
&=& \left (\frac{\beta_3}{4}-  \frac{1}{4}+\frac{ \alpha_0}{2} -\frac{\alpha_3}{4} \right) x \bar{L}_{3n} +(\alpha_0-1) \left (\frac{\beta_3}{4}-  \frac{1}{4}+\frac{ \alpha_0}{2} -\frac{\alpha_3}{4} \right)  x^2 \bar{L}^2_{3n}  \nonumber\\
&& + \frac{O(1)\Delta^{3,5}_{n,x}}{x^2}.
\end{eqnarray*}
For the other terms in (\ref{CL46}) involving $x^2 \bar{L}_{4n}$, $x^2 \bar{L}^2_{3n}$ and $x^2B_n^{-4}\sum_{i=1}^n ( E \bar{X}_i^2)^2$, by Lemma \ref{properties}(i), $\lambda_{\delta}$ can be replaced by $1/2$ with a difference of $O(1)\Delta^{3,5}_{n,x}/x^2$. Therefore (\ref{CL42}) follows. 
\end{proof}

\begin{proposition}\label {corollary1}
Let $S_{n, (k)} =\sum_{i=1}^n X_{i, (k)}$ and $V_{n, (k)}^2= \sum_{i=1}^n X^2_{i,(k)}$, where $X_{i, (1)}=X_i$ and $X_{i, (2)} =X_i I(|bX_i| \leq 1)$. Under conditions (\ref{condition1}) and (\ref{condition2}),  for $x > 0$, $|\delta(x)| \leq x^2/5$ and $k=1$ or $2$, 
\begin{eqnarray}\label {DG83a}
&& P \bigg \{ 2bS_{n,(k)} - (\alpha_0- \alpha_3 x \bar{L}_{3n}) b^2 V^2_{n,(k)}   \geq (2-\alpha_0)x^2 +\delta(x) \bigg \} \nonumber\\
&=&  \exp \Bigg \{\frac{x^2}{2}+\sum_{i=1}^n \log E e^{\lambda_{\delta} \xi( X_{i, (k)})}- \lambda_{\delta} (2-\alpha_0)x^2 - \lambda_{\delta}\delta(x)    \Bigg \} (1-\Phi(x))\nonumber\\
&&\times  \left(1+O(1) x\mathcal{L}_{3n}\right).
\end{eqnarray}
If $\delta(x)=\beta_3 x^3 \bar{L}_{3n}+\beta_4 x^4 \bar{L}_{4n}+ \beta_5 x^4 \bar{L}^2_{3n}+O( 1)\Delta^{3,5}_{n,x}$ with $|\beta_j|\leq \sqrt{A}$, then
\begin{eqnarray}\label {22AH}
&& \frac{x^2}{2}+\sum_{i=1}^n \log E e^{\lambda_{\delta} \xi( X_{i, (k)})}- \lambda_{\delta} (2-\alpha_0)x^2 - \lambda_{\delta}\delta(x) \nonumber\\
&=&   \left (\frac{\alpha_3-\beta_3}{2}+\frac{1-3\alpha_0}{6} \right )x^3 \bar{L}_{3n} + \left (-\frac{\beta_4}{2} +\frac{(\alpha_0-1)^2}{8}- \frac{1}{12} \right )x^4 \bar{L}_{4n}  \nonumber\\
&& + \bigg ( \frac{\alpha_3-\beta_5}{2} - \frac{\big (\beta_3 -\alpha_3+ 2 \alpha_0 -1\big )^2}{8} \bigg )  x^4 \bar{L}^2_{3n}  \nonumber\\
&& -\frac{\left(1-  \alpha_0 \right)^2}{8} \frac{x^4\sum_{i=1}^n ( E X_i^2)^2}{B_n^4} +O (1) \Delta^{3,5}_{n,x}  .
\end{eqnarray}
\end{proposition}
\begin{proof}
By the conjugate method (e.g., (4.9) of Petrov \cite{Petrov65} or (2.11) of Petrov \cite{Petrov}) together with (\ref{CD30}) and (\ref{CG39}), 
\begin{eqnarray} \label {50al}
&& P \bigg \{ 2bS_{n,(k)} - (\alpha_0+ \alpha_3 x \bar{L}_{3n}) b^2 V^2_{n,(k)}  \geq (2-\alpha_0)x^2 +\delta(x) \bigg \} \nonumber\\
&=& \int \cdots \int   I  \bigg \{ \sum_{i=1}^n \xi(x_i) \geq (2-\alpha_0)x^2  + \delta(x)  \bigg  \} \prod_{i=1}^n d P( X_{i, (k)} \leq x_i) \nonumber\\
&=& \int \cdots \int  e^{ {\sum_i \log E e^{\lambda_{\delta} \xi (X_{i,(k)})}}-\lambda_{\delta}\sum_i \xi(x_i) } \nonumber\\
&&  \qquad \quad  \times I  \bigg \{ \sum_{i=1}^n \xi(x_i)  \geq (2-\alpha_0)x^2  + \delta(x)  \bigg \} \prod_{i=1}^n d P(Z_{\lambda_{\delta},i, (k)} \leq x_i) \nonumber\\
&=& e^{ {\sum_i \log E e^{\lambda_{\delta} \xi ( X_{i, (k)})}}-\lambda_{\delta} (2 -\alpha_0)x^2 - \lambda_{\delta} \delta(x)  }  E\bigg \{ e^{-\lambda_{\delta} \sum_i \{ \xi (Z_{\lambda_{\delta}, i, (k)}) - E \xi (Z_{\lambda_{\delta}, i, (k)})\}} \nonumber\\
&& \times  I \Big ( \sum_{i=1}^n ( \xi (Z_{\lambda_{\delta}, i, (k)}) - E  \xi (Z_{\lambda_{\delta}, i, (k)})) \geq 0 \Big )\bigg \}. \qquad \ 
\end{eqnarray}
Let $W_i =\xi (Z_{\lambda_{\delta}, i, (k)})$ and $S = \sum_{i=1}^n (W_i - EW_i)$. By (\ref{Vareta1}), we have 
$$B_n'^2:= \sum_{i=1}^nE (W_i-EW_i)^2  =\sum_{i=1}^n Var(\xi (Z_{\lambda_{\delta}, i, (k)})) = 4x^2 +O(1)x^3 \mathcal{L}_{3n}. $$ Then by the proof of Proposition 2.2 of Wang \cite{Wang} with $\theta=0$, 
\begin{eqnarray}\label {22JK}
&& E\bigg \{ e^{-\lambda_{\delta} \sum \{ \xi (Z_{\lambda_{\delta}, i, (k)}) - E \xi (Z_{\lambda_{\delta}, i, (k)})\}} \nonumber\\
&& \times  I \Big ( \sum_{i=1}^n ( \xi (Z_{\lambda_{\delta}, i, (k)}) - E  \xi (Z_{\lambda_{\delta}, i, (k)})) \geq 0 \Big )\bigg \}\nonumber\\
&=& E e^{-\lambda_{\delta}S}I(S \geq 0) \nonumber\\
& =&   e^{x^2/2} (1-\Phi(x)) (1+O(1)x\mathcal{L}_{3n}) .
\end{eqnarray}
Hence (\ref{DG83a}) follows from (\ref{50al}) and (\ref{22JK}). By (\ref{noxi}), 
\begin{eqnarray} \label{DG85}
&&\sum_{i=1}^n \log E e^{\lambda_{\delta} \xi(X_{i, (k)})}- \lambda_{\delta} (2-\alpha_0)x^2\nonumber\\
&=&(2 \lambda_{\delta} ^2-2\lambda_{\delta}  )x^2  + \left (\frac{4  \lambda_{\delta} ^3}{3} -2 \lambda_{\delta} ^2 \alpha_0 +\lambda_{\delta}  \alpha_3 \right )x^3 \bar{L}_{3n} \nonumber\\
&& + \left ( \frac{2\lambda_{\delta} ^4}{3}-2 \lambda_{\delta} ^3 \alpha_0+\frac{\lambda_{\delta} ^2 \alpha_0^2}{2} \right )x^4 \bar{L}_{4n} + 2 \lambda_{\delta} ^2 \alpha_3 x^4 \bar{L}^2_{3n}  \nonumber\\
&&-\frac{\left(2 \lambda_{\delta} ^2-\lambda_{\delta}  \alpha_0 \right)^2}{2} \frac{x^4\sum_{i=1}^n ( E \bar{X}_i^2)^2}{B_n^4} +O(1) \Delta^{3,5}_{n,x}.
\end{eqnarray}
From the expression of $\lambda_{\delta}$ in (\ref{CL42}) and by Lemma \ref{properties}(i), we have
\begin{eqnarray}\label {DG86}
(2 \lambda_{\delta} ^2-2\lambda_{\delta}  )x^2 = -\frac{x^2}{2}+ 2 \left (\frac{\beta_3}{4}+ \frac{\alpha_0}{2} - \frac{\alpha_3}{4} -\frac{1}{4}\right )^2     x^4 \bar{L}^2_{3n} +O(1)\Delta^{3,5}_{n,x} 
\end{eqnarray}
and
\begin{eqnarray}\label {DG87a}
&& \left (\frac{4  \lambda_{\delta} ^3}{3} -2 \lambda_{\delta} ^2 \alpha_0 +\lambda_{\delta}  \alpha_3 \right )x^3 \bar{L}_{3n}\nonumber\\
&=& \left (\frac{1}{6}-\frac{\alpha_0}{2} +\frac{\alpha_3}{2} \right )x^3 \bar{L}_{3n}+ \left (  1-2 \alpha_0 +\alpha_3 \right )\left (\frac{\beta_3}{4}+ \frac{\alpha_0}{2} - \frac{\alpha_3}{4} -\frac{1}{4}\right )     x^4 \bar{L}^2_{3n} \nonumber\\
&& + O(1)\Delta^{3,5}_{n,x}.
\end{eqnarray}
Applying (\ref{DG86}) and (\ref{DG87a}) to (\ref{DG85}) and by Lemma \ref{properties}(i), 
\begin{eqnarray} \label {DG87}
&&\sum_{i=1}^n \log E e^{\lambda_{\delta} \xi(X_{i, (k)})}- \lambda_{\delta} (2-\alpha_0)x^2\nonumber\\
&=& -\frac{x^2}{2} + \left (\frac{1}{6}-\frac{\alpha_0}{2} +\frac{\alpha_3}{2} \right )x^3 \bar{L}_{3n} + \left ( \frac{1}{24}- \frac{\alpha_0}{4}+\frac{\alpha_0^2}{8} \right )x^4 \bar{L}_{4n}  \nonumber\\
&&+ 2 \left (\frac{\beta_3}{4}+ \frac{\alpha_0}{2} - \frac{\alpha_3}{4} -\frac{1}{4}\right )^2     x^4 \bar{L}^2_{3n}\nonumber\\
&&+ \left (  1-2 \alpha_0 +\alpha_3 \right )\left (\frac{\beta_3}{4}+ \frac{\alpha_0}{2} - \frac{\alpha_3}{4} -\frac{1}{4}\right )     x^4 \bar{L}^2_{3n} + \frac{\alpha_3}{2}x^4 \bar{L}^2_{3n} \nonumber\\
&&-\frac{\left(1-  \alpha_0 \right)^2}{8} \frac{x^4\sum_{i=1}^n ( E \bar{X}_i^2)^2}{B_n^4} +O(1) \Delta^{3,5}_{n,x}.
\end{eqnarray}
Since $\delta(x)=\beta_3 x^3 \bar{L}_{3n}+\beta_4 x^4 \bar{L}_{4n}+ \beta_5 x^4 \bar{L}^2_{3n}+O(1)\Delta^{3,5}_{n,x}$, then
\begin{eqnarray} \label {DG88}
-\lambda_{\delta}\delta(x) &=&- \frac{\beta_3}{2} x^3 \bar{L}_{3n}-\frac{\beta_4}{2} x^4 \bar{L}_{4n}\nonumber\\
&& - \left \{\frac{\beta_5}{2} +\beta_3\left (\frac{\beta_3}{4}+ \frac{\alpha_0}{2} - \frac{\alpha_3}{4} -\frac{1}{4}\right )   \right \}x^4 \bar{L}^2_{3n}+O(1)\Delta^{3,5}_{n,x}. \qquad
\end{eqnarray}
Combining (\ref{DG87}) and (\ref{DG88}), we obtain (\ref{22AH}). 
\end{proof}


\begin{proposition} \label {Proposition2}
Under conditions (\ref{condition1}) and (\ref{condition2}), for $x \geq 2$,
\begin{eqnarray} \label {AQ65}
P (S_{n} \geq x V_{n}) &\leq& A \exp \bigg \{-\frac{1}{3}x^3\bar{L}_{3n} - \frac{1}{12}  x^4\bar{L}_{4n}+A\Delta^{3,5}_{n,x} \bigg \} \nonumber\\
&& \times (1-\Phi(x))(1+Ax \mathcal{L}_{3n}).
\end{eqnarray}
Consequently, for $S_n^{(i)}=S_n -X_i$ and $V_n^{(i)}= (V_n^2-X^2_i)^{1/2}$,
\begin{eqnarray} \label {AR68}
P \left (S_n^{(i)} \geq (x^2-1)^{1/2} V^{(i)}_n \right) &\leq& A\exp \bigg \{-\frac{1}{3}x^3\bar{L}_{3n} - \frac{1}{12}  x^4\bar{L}_{4n}+A\Delta^{3,5}_{n,x} \bigg \} \nonumber\\
&&\times (1-\Phi(x))(1+Ax \mathcal{L}_{3n}).
\end{eqnarray}
 \end{proposition}
\begin{proof}
Let $P (S_n \geq x V_n) = I_1+I_2+I_3+I_4$ where
\begin{eqnarray}
I_1  &=& P \left \{S_n \geq x V_n, \ \     | V_n/B_n -1|  \leq  |x \bar{L}_{3n}|/2+x^2 \bar{L}_{4n} \right \}, \nonumber \\
I_2  &=& P \left \{S_n \geq x V_n, \ \   0 \leq  V_n/B_n   \leq 1- |x \bar{L}_{3n}|/2-x^2\bar{L}_{4n} \right \}, \nonumber \\
I_3  &=& P \left \{S_n \geq x V_n, \ \   1+ |x \bar{L}_{3n}|/2+x^2\bar{L}_{4n} < V_n/B_n \leq 3 \right \}, \nonumber \\
I_4 &=& P \left \{S_n \geq x V_n, \ \   V_n/B_n > 3 \right \}. \nonumber
\end{eqnarray}
Observe that for any real valued set $B$ and any constants $\gamma>0$ and $\theta>0$,
\begin{eqnarray}\label {BT70}
&& P \left \{ S_n \geq x V_n, \ \frac{V_n}{B_n} \in  B \right \} \nonumber\\
&\leq& P \left \{\gamma bS_n -\theta b^2 V_n^2\geq \inf_{V_n/B_n \in B} \left (\gamma x \frac{xV_n}{B_n} - \theta  \frac{x^2 V_n^2}{B_n^2} \right)\right \}.
\end{eqnarray}
Let $t=xV_n / B_n$. Then $f_x (t):=2 x t - t^2$ has minimum at $a$ for $t \in [a, x]$ and has minimum at $b$ for $t \in [x, b]$. Hence by (\ref{BT70}), Lemma \ref{properties}(i)  and Proposition \ref{corollary1} with $\alpha_0=1, \ \alpha_3 = 0, \ \beta_3=0, \ \beta_4=0$ and $\beta_5 = -1/4$,
\begin{eqnarray}
 I_1\leq P \bigg \{2bS_n - b^2 V_n^2\geq  2 x^2 \bigg (1-\frac{|x \bar{L}_{3n}|}{2}-x^2\bar{L}_{4n}  \bigg ) -  x^2 \bigg (1-\frac{|x \bar{L}_{3n}|}{2}- x^2\bar{L}_{4n} \bigg ) ^2 \bigg \}\nonumber\\
+P \bigg \{2bS_n - b^2 V_n^2\geq  2 x^2 \bigg (1+\frac{|x \bar{L}_{3n}|}{2}+x^2\bar{L}_{4n}  \bigg ) -  x^2 \bigg (1+\frac{|x \bar{L}_{3n}|}{2}+ x^2\bar{L}_{4n} \bigg ) ^2 \bigg \}\nonumber\\
= 2 P \bigg \{2bS_n - b^2 V_n^2\geq   x^2 - \frac{1}{4}x^4\bar{L}^2_{3n}+O(1)\Delta^{3,5}_{n,x} \bigg \} \quad \qquad  \qquad \qquad  \qquad \qquad  \qquad \ \  \nonumber\\
\leq 2\exp \bigg \{-\frac{1}{3}x^3\bar{L}_{3n}-\frac{1}{12} x^4\bar{L}_{4n} +A\Delta^{3,5}_{n,x}\bigg \} (1-\Phi(x))(1+Ax \mathcal{L}_{3n}).   \qquad \qquad  \qquad  \nonumber
\end{eqnarray}

To estimate $I_2$, note that $f_x(t):=2xt-7t^2/3$ has minimum at $a$ for $t \in [0, a]$ with $a \geq 6x/7$. Then by (\ref{BT70}), Lemma \ref{properties}(i) and Proposition \ref{corollary1} with $\alpha_0=7/3, \ \alpha_3 = 0, \ \beta_4=8/3, \ \beta_5 = -7/12$, and $\beta_3=4/3$ if $x^3 \bar{L}_{3n}\geq 0$ and $\beta_3=-4/3$ if $x^3 \bar{L}_{3n}< 0$, 
\begin{eqnarray}
I_2 &\leq&  P \bigg \{2bS_n - \frac{7}{3} b^2 V_n^2\geq  2 x^2 \bigg (1-\frac{|x \bar{L}_{3n}|}{2}-x^2\bar{L}_{4n}  \bigg ) \nonumber\\
&&-  \frac{7}{3}x^2 \bigg (1-\frac{|x \bar{L}_{3n}|}{2}- x^2\bar{L}_{4n} \bigg ) ^2 \bigg \}\nonumber\\
&=&  P \bigg \{2b S_n - \frac{7}{3} b^2 V_n^2\geq  -\frac{ x^2 }{3}+\frac{4}{3} |x^3\bar{L}_{3n}| +\frac{8}{3}x^4 \bar{L}_{4n}-\frac{7}{12}x^4 \bar{L}^2_{3n}+O(1)\Delta^{3,5}_{n,x}\bigg \}\nonumber\\
&\leq& \exp \bigg \{-\frac{1}{3}x^3\bar{L}_{3n} - \frac{1}{12}  x^4\bar{L}_{4n}+A\Delta^{3,5}_{n,x} \bigg \}(1-\Phi(x))(1+Ax \mathcal{L}_{3n}) \nonumber
\end{eqnarray}

To estimate $I_3$, note that $f_x(t):=2xt-t^2/3$ has minimum at $a$ for $t \in [a, 3x]$. Then by (\ref{BT70}), Lemma \ref{properties}(i) and  Proposition \ref{corollary1} with $\alpha_0=1/3, \ \alpha_3=0, \ \beta_4=4/3, \ \beta_5=-1/12$, and $\beta_3=2/3$ if $x^3 \bar{L}_{3n}\geq 0$ and $\beta_3=-2/3$ if $x^3 \bar{L}_{3n}< 0$, 
\begin{eqnarray}
I_3 &\leq&  P \bigg \{2 b S_n - \frac{1}{3} b^2 V_n^2\geq  2 x^2 \bigg (1+\frac{|x \bar{L}_{3n}|}{2} +x^2\bar{L}_{4n} \bigg ) \nonumber\\
&&-  \frac{1}{3}x^2 \bigg (1+\frac{|x \bar{L}_{3n}|}{2} +x^2 \bar{L}_{4n}\bigg ) ^2 \bigg \}\nonumber\\
&=&  P \bigg \{2 b S_n - \frac{1}{3} b^2 V_n^2\geq  \frac{5 x^2 }{3}+\frac{2}{3} |x^3\bar{L}_{3n}|+\frac{4}{3}x^4 \bar{L}_{4n} -\frac{1}{12}x^4 \bar{L}^2_{3n}+O(1) \Delta^{3,5}_{n,x}\bigg \} \nonumber \\
&\leq& \exp \bigg \{-\frac{1}{3}x^3\bar{L}_{3n} - \frac{1}{12}  x^4\bar{L}_{4n}+A\Delta^{3,5}_{n,x} \bigg \} (1-\Phi(x))(1+Ax \mathcal{L}_{3n}), \nonumber
\end{eqnarray}
where we also use the fact that $\bar{L}^2_{3n}\leq \bar{L}_{4n}$. 

Finally by the proof of Lemma 8.1 of Jing, Shao and Wang \cite{JingShaoWang}, for $x \geq 2$,
\begin{eqnarray}
I_4 \leq 2 e^{-x^2}\leq A\exp \bigg \{-\frac{1}{3}x^3\bar{L}_{3n} - \frac{1}{12}  x^4\bar{L}_{4n}+A\Delta^{3,5}_{n,x} \bigg \} (1-\Phi(x)).\nonumber
\end{eqnarray}
Therefore, (\ref{AQ65}) follows from all of the above. 
\end{proof}
\bigskip


For $x \geq 2$, let $\tau=B_n/x$ and set
\begin{gather}
\bar{X}_i=X_i I(|X_i| \leq \tau), \quad  \bar{Z}_{\lambda_{\delta},i} = Z_{\lambda_{\delta}, i, (2)}, \quad \bar{S}_n = \sum_{i=1}^n \bar{X}_i, \quad  \bar{V}_n^2 = \sum_{i=1}^n \bar{X}_i^2,  \nonumber\\
  \bar{V}_{n, \delta}^{2} = \sum_{i=1}^n  \bar{Z}_{\lambda_{\delta}, i}^2, \quad u_{i}= \bar{Z}_{\lambda_{\delta},i}^2-E  \bar{Z}_{\lambda_{\delta},i}^2, \quad \bigtriangledown_n=\frac{x^2}{B_n^4}\Big (\sum_{i=1}^n u_i \Big)^2.    \label {81tg}
\end{gather}


\begin{proposition}\label {Prop2.2}
Under conditions (\ref{condition1}) and (\ref{condition2}), for $x \geq 2$,
\begin{eqnarray}
P (\bar{S}_n \geq x \bar{V}_n) \leq \exp \Big \{-\frac{1}{3}x^3 \bar{L}_{3n}  -  \frac{1}{12}  x^4\bar{L}_{4n}+A \Delta^{3,5}_{n,x}\Big \} (1-\Phi(x))(1+A x \mathcal{L}_{3n}).\nonumber
\end{eqnarray}
\end{proposition}
\begin{proof}
Let $ \xi(y)=2b y - (1-x\bar{L}_{3n}/2)b^2 y^2$ and
\begin{eqnarray}\label {delta(x)}
\delta(x) = \frac{x^3 \bar{L}_{3n}}{2} + \frac{x^4 \bar{L}^2_{3n} }{4}-2A\Delta^{3,5}_{n,x}.
\end{eqnarray}
By the conjugate method similar to (\ref{50al}),
\begin{eqnarray}\label {66tj}
&& P (\bar{S}_n \geq x \bar{V}_n)\nonumber\\
&=& \int \cdots \int   I  \bigg \{ \sum_{i=1}^n x_i \geq x \Big (\sum_{i=1}^n x_i^2\Big)^{1/2} \bigg  \} \prod_{i=1}^n d P( \bar{X}_{i} \leq x_i) \nonumber\\
&=& \int \cdots \int  e^{ {\sum_i \log E e^{\lambda_{\delta} \xi (\bar{X}_i)}}-\lambda_{\delta}\sum_i \xi(x_i) } I  ( s_n \geq x v_n ) \prod_{i=1}^n d P(\bar{Z}_{\lambda_{\delta},i} \leq x_i) , \qquad
\end{eqnarray}
where $s_n = \sum_{i=1}^n x_i $ and $v^2_n = \sum_{i=1}^n x_i^2$. We can write
\begin{eqnarray}\label {68tj}
v_n = \frac{ B_n}{1- x \bar{L}_{3n}/2} \left (1+ y_n\right )^{1/2},
\end{eqnarray}
where
\begin{eqnarray}\label {47tp}
y_n = \frac{(1- x \bar{L}_{3n}/2)^2 v_n^2-B_n^2 }{B_n^2} . 
\end{eqnarray}
Since $\bar{L}^2_{3n} \leq \bar{L}_{4n}$, then $\lambda_{\delta}= 1/2 + O(1) x \bar{L}_{3n}+O(1)x^{-2}\Delta^{3,5}_{n,x}$ by (\ref{CL42}. Hence by (\ref{CD31}) and (\ref{81tg}),  
\begin{eqnarray} \label {74tj}
E \bar{V}_{n,\delta}^{2}&=& \sum_{i=1}^n E\bar{X}_i^2 +\frac{x}{B_n}\sum_{i=1}^n E \bar{X}_i^3 +O(1)\frac{x^2}{B_n^2}\sum_{i=1}^n E \bar{X}_i^4 +O(1)x^{-2}\Delta^{3,5}_{n,x}B_n^2 \nonumber\\
&=& B_n^2 + x \bar{L}_{3n}B_n^2 + O(1) x^2 \bar{L}_{4n} B_n^2 + O(1) x^{-2}\Delta^{3,5}_{n,x}B_n^2. \nonumber
\end{eqnarray}
Since $(1-  x \bar{L}_{3n}/2)^2= 1-  x \bar{L}_{3n} +  x^2 \bar{L}^2_{3n}/4 $. Then
\begin{eqnarray}\label {49tp}
  (1-  x \bar{L}_{3n}/2)^2 E  \bar{V}_{n,\delta}^{2}  = B_n^2 +O(1) x^2 \bar{L}_{4n} B_n^2 + O(1) x^{-2}\Delta^{3,5}_{n,x}B_n^2.
\end{eqnarray}
Applying (\ref{49tp}) to (\ref{47tp}), we have
\begin{eqnarray}\label {50tp}
y_n &=& \frac{(1- x \bar{L}_{3n}/2)^2 (v_n^2-E  \bar{V}^2_{n,\delta}) }{B_n^2}+\frac{(1- x \bar{L}_{3n}/2)^2 E  \bar{V}_{n,\delta}^2-B_n^2}{B_n^2}\nonumber\\
&=&\frac{(1- x \bar{L}_{3n}/2)^2 (v_n^2-E  \bar{V}^2_{n,\delta}) }{B_n^2}+O(1)x^2 \bar{L}_{4n} +O(1)x^{-2}\Delta^{3,5}_{n,x} .
\end{eqnarray}
Note that $x  \bar{L}^2_{4n} \leq x \bar{\mathcal{L}}_{3n }\mathcal{\bar{L}}_{5n} \leq \mathcal{\bar{L}}_{5n}/A$ by Lemma \ref{properties}(i). Then
\begin{eqnarray}\label {51tp}
y_n^2 &\leq&  \frac{2(1- x \bar{L}_{3n}/2)^4  (v_n^2-E  \bar{V}_{n,\delta}^2)^2 }{B_n^4}+ Ax^{-2}\Delta^{3,5}_{n,x}.
\end{eqnarray}
Observe that $(1+y)^{1/2} \geq 1+y/2-y^2/2$ for any $y \geq -1$ and $(1+y)^{1/2} \geq 1+y/2-y^2/m$ for $y \geq m/2$ with $m >0$. By (\ref{50tp}), $ B_n^{-2}(v_n^2-E  \bar{V}^2_{n,\delta}) \geq -3/2$ because $y_n \geq -1$. Moreover, $B_n^{-2}(v_n^2-E  \bar{V}_{n,\delta}^2) \geq m$ implies that $y_n \geq m/2$ for $m \geq 2$. Hence
\begin{eqnarray}\label {52tp}
(1+y_n)^{1/2} &\geq& 1 + \frac{y_n}{2}- \theta  y_n^2,
\end{eqnarray}
where
\[ \theta =
  \begin{cases}
    1/2       & \quad \text{if } B_n^{-2}|v_n^2-E  \bar{V}^2_{n,\delta}| \leq 2  \\
    1/m  & \quad \text{if } m < B_n^{-2}(v_n^2-E  \bar{V}^2_{n,\delta}) \leq m+1 \ \textrm{for $m \geq 2$}. \\
  \end{cases}
\]
Since $4(1- x \bar{L}_{3n}/2)^3 \leq 5$, then combining (\ref{68tj}), (\ref{47tp}), (\ref{51tp}) and (\ref{52tp}), we have
\begin{eqnarray}\label {53tp}
 && I  ( s_n \geq x v_n   ) =  I \Big \{  2 bs_n \geq \frac{2x^2}{1- x \bar{L}_{3n}/2} \left (1+ y_n\right )^{1/2}    \Big \}\nonumber\\
&\leq& I \bigg \{ 2 b s_n - (1- x \bar{L}_{3n}/2)b^2v_n^2 + \frac{5 \theta x^2(v_n^2-E  \bar{V}^2_{n,\delta})^2 }{B_n^4}  \nonumber\\
&& \geq  \frac{x^2}{1- x \bar{L}_{3n}/2}- A \Delta^{3,5}_{n,x}\bigg \}. 
\end{eqnarray}
From (\ref{DC109}) and (\ref{delta(x)}),
\begin{eqnarray}\label {78tja}
\frac{ x^2}{1- x \bar{L}_{3n}/2}-   A \Delta^{3,5}_{n,x} \geq x^2 + \delta(x) . 
\end{eqnarray}
By (\ref{66tj}), (\ref{53tp}), (\ref{78tja}) and the definition of $\bigtriangledown_n$ in (\ref{81tg}),
\begin{eqnarray} 
&& P (\bar{S}_n \geq x \bar{V}_n)\nonumber\\
&\leq& E    e^{ {\sum_i \log E e^{\lambda_{\delta} \xi (\bar{X}_i)}}-\lambda_{\delta}\sum_i \xi(\bar{Z}_{\lambda_{\delta}, i}) } I  \Big \{\sum_{i=1}^n \xi (\bar{Z}_{\lambda_{\delta}, i})+5 \theta \bigtriangledown_{n} \geq  x^2+\delta(x) \Big \}  . \qquad 
\end{eqnarray}
Let
\begin{eqnarray}
T_{n } = \sum_{i=1}^n (\xi ( \bar{Z}_{\lambda_{\delta}, i })- E \xi (\bar{Z}_{\lambda_{\delta}, i }))+5 \theta \bigtriangledown_{n}.\nonumber
\end{eqnarray}
Note that $\sum_{i=1}^n E \xi (\bar{Z}_{\lambda_{\delta}, i })= x^2 +\delta(x)$ in (\ref{CG39}). Hence
\begin{eqnarray}\label {82tj}
 P (\bar{S}_n \geq x \bar{V}_n)&\leq&  \exp \Big \{ {\sum \log E e^{\lambda_{\delta} \xi (\bar{X}_i)}}-\lambda_{\delta}x^2 -\lambda_{\delta} \delta(x) \Big\} \nonumber\\
&& \times E   e^{-\lambda_{\delta} \sum_i (\xi( \bar{Z}_{\lambda_{\delta}, i})- E \xi ( \bar{Z}_{\lambda_{\delta}, i}))}  I \left  ( T_{n }  \geq 0 \right  ) .
\end{eqnarray}
Since $|e^s - 1 | \leq e^{0 \vee s}|s|$, then
\begin{eqnarray}\label {J123}
&& E   e^{-\lambda_{\delta} \sum_i (\xi( \bar{Z}_{\lambda_{\delta}, i})- E \xi ( \bar{Z}_{\lambda_{\delta}, i}))}  I \left  ( T_{n }  \geq 0 \right  ) \nonumber\\
&=& E   e^{ 5 \theta \lambda_{\delta}\bigtriangledown_{n }-\lambda_{\delta}T_{n } }  I \left  ( T_{n }  \geq 0 \right  ) \nonumber\\
&\leq&E   e^{ -\lambda_{\delta}T_{n } }  I \left  ( T_{n }  \geq 0 \right  ) +5 \theta \lambda_{\delta} E \bigtriangledown_{n } e^{5 \theta  \lambda_{\delta}\bigtriangledown_{n }-\lambda_{\delta}T_{n} }  I \left  ( T_{n }  \geq 0 \right  ) \nonumber\\
&\leq&E   e^{ -\lambda_{\delta}T_{n } }  I \left  ( T_{n }  \geq 0 \right  ) +5 \theta \lambda_{\delta} E\bigtriangledown_{n } e^{ 5 \theta \lambda_{\delta}\bigtriangledown_{n }}  \nonumber\\
&:=& J_1 +J_2. 
\end{eqnarray}
To estimate $J_1$, let $B_n'^2=\sum_{i =1}^n (K_i-E K_i)^2$ where $K_i =  \xi (\bar{Z}_{\lambda_{\delta},i}) + 5 \theta x^2  B_n^{-4}u_i^2$. Then $B_n'^2 = 4x^2 +O(x^3 \mathcal{L}_{3n})$ by (\ref{Vareta1}) and (\ref{CD31}). Following the proof of Proposition 2.2 of Wang \cite{Wang}, we have
\begin{eqnarray}\label {J1}
J_1 =e^{x^2/2} (1-\Phi(x)) (1+ A x \mathcal{L}_{3n}).
\end{eqnarray}
By Lemma \ref{Lemma3.4} below and Lemma \ref{properties}(ii),
\begin{eqnarray}\label {94ke}
J_2 \leq  A x^2 \bar{\mathcal{L}}_{4n} \leq  A ( \mathcal{\bar{L}}_{3n}+ x^{-1}\Delta^{3,5}_{n,x}). 
\end{eqnarray}
Note that $(2 \pi)^{-1/2}(x^{-1}-x^{-3})e^{-x^2/2} \leq (1-\Phi(x))$ for $x \geq 2$. Then by (\ref{J123})-(\ref{94ke}),
\begin{eqnarray}\label {62ah}
&& E  e^{-\lambda_{\delta} \sum_i \xi(Z_{\lambda_{\delta}, i, (2)})- E \xi (Z_{\lambda_{\delta}, i, (2)}))}  I \left  ( T_n  \geq 0 \right  )   \nonumber\\
&\leq& e^{x^2/2} (1-\Phi(x)) (1+Ax \mathcal{L}_{3n}+ A \Delta^{3,5}_{n,x}) \nonumber\\
&\leq& e^{x^2/2+A\Delta^{3,5}_{n,x}} (1-\Phi(x)) (1+Ax \mathcal{L}_{3n}).
\end{eqnarray}
By (\ref{22AH}) in Proposition \ref{corollary1} and the $\delta(x)$ in (\ref{delta(x)}), 
\begin{eqnarray}\label {67EL}
 \frac{x^2}{2}+\sum_{i=1}^n \log E e^{\lambda_{\delta} \xi(\bar{X}_i)}- \lambda_{\delta} x^2 - \lambda_{\delta}\delta(x) \leq -\frac{1}{3} x^3 \bar{L}_{3n} -\frac{1}{12}x^4 \bar{L}_{4n} + A\Delta^{3,5}_{n,x}. \ \ 
\end{eqnarray}
Therefore, the proposition follows from (\ref{82tj}), (\ref{62ah}) and (\ref{67EL}).
\end{proof}


\begin{lemma}\label {Lemma3.4}
Suppose that $\theta = 1/2$ if $B_n^{-2}|\sum_{i=1}^n u_i| \leq 2$ and $\theta = 1/m$ if $m < B_n^{-2}(\sum_{i=1}^n u_i) \leq m+1$ for $m \geq 2$. Then for $x \geq 2$ and  $\bigtriangledown_{n } $ defined in (\ref{81tg}),
\begin{eqnarray}
E \bigtriangledown_{n } e^{5 \theta   \bigtriangledown_{n}}  \leq A x^2 \bar{\mathcal{L}}_{4n}. \nonumber
\end{eqnarray}
\end{lemma}
\begin{proof}
By the definition of $\bigtriangledown_{n } $ in (\ref{81tg}),
\begin{eqnarray} \label {44th}
 E \bigtriangledown_{n }   e^{5 \theta   \bigtriangledown_{n }}=\frac{x^2}{B_n^4} \sum_{i=1}^n E u_i^2 e^{5 \theta \bigtriangledown_{n}}   + \frac{x^2}{B_n^4} \sum_{1 \leq i, j \leq n, i \neq j} E u_i u_j  e^{5 \theta \bigtriangledown_{n}}.
\end{eqnarray}
First, we want to show that for $d=\{1,...,n \}\backslash \{i\}$ or $d=\{1,...,n \}\backslash \{i, j\}$,
\begin{eqnarray}\label {66tq}
E e^{10 \theta x^2 B_n^{-4}(\sum_{k \in d}u_k)^2} \leq A.
\end{eqnarray}
Since $B_n^{-2}|u_k| \leq 1/x^2$ and since $ B_n^{-2}(\sum_{k=1}^n u_k) \geq -3/2$ as mentioned after (\ref{51tp}), then $ B_n^{-2}(\sum_{k\in d}^n u_k) \geq -3/2 -2/x^2 \geq -2$. Hence
\begin{eqnarray}\label {67tq}
&& E e^{10 \theta x^2 B_n^{-4} (\sum_{k \in d}u_k)^2} \nonumber\\
&=& E e^{10 \theta x^2 B_n^{-4}(\sum_{k \in d}u_k)^2} I\Big (x B_n^{-2} \Big|\sum_{k\in d}  u_k \Big| \leq 2 \Big) \nonumber\\
&&+ E e^{10 \theta x^2 B_n^{-4}(\sum_{k \in d}u_k)^2} I\Big (2 < x B_n^{-2}\Big |\sum_{k\in d}  u_k\Big |\leq 2 x \Big)  \nonumber\\
&&  + \sum_{m=2}^{\infty}E e^{10 \theta x^2 B_n^{-4}(\sum_{k \in d}u_k)^2} I\Big (m x < x B_n^{-2} \sum_{k \in d} u_k \leq (m+1)x \Big) \nonumber\\
&:=& K_1+K_2 +K_3. 
\end{eqnarray}
Since $\theta \leq 1/2$, then
\begin{eqnarray} \label {68tq}
K_1 \leq e^{20}. 
\end{eqnarray}
Let $[2x]$ be the largest integer smaller than or equal to $x$. Then
\begin{eqnarray}\label {68tqa}
K_2 &\leq& \sum_{l =2}^{[2x]} E e^{10 \theta x^2 B_n^{-4}(\sum_{k \in d}u_k)^2} I\Big (l  <  x B_n^{-2}\Big |\sum_{k\in d} u_k \Big | \leq l+1 \Big) \nonumber\\
&\leq&  \sum_{l =2}^{[2x]} e^{5  (l+1)^2} P\Big (  x B_n^{-2}\Big |\sum_{k \in d }  u_k \Big |> l \Big). 
\end{eqnarray}
Since $|e^{s}-1-s| \leq e^{0 \vee s}s^2/2$ and $E u_k=0$, then for $t>0$, 
\begin{eqnarray}
P\Big (  x B_n^{-2}\Big |\sum_{k \in d}  u_k \Big |> l \Big) &\leq&  e^{-tl} \prod_{k\in d}  E e^{t x B_n^{-2}  u_k}+e^{-tl} \prod_{k\in d}  E e^{-t x B_n^{-2}  u_k}  \nonumber\\
&\leq&   2 e^{-tl} \prod_{k \in d} \Big(1+\frac{1}{2} e^{tx B_n^{-2}|u_k|} t^2x^2 B_n^{-4}E u_k^2 \Big) .  \nonumber
\end{eqnarray}
Let $t= (\log A^{1/2})l/2$. Since $x B_n^{-2}|u_k| \leq 1/x$, then $e^{t xB_n^{-2}|u_k|} \leq A^{1/2}$ for $l \leq 2 x$. Note that $E u_k^2 \leq E \bar{Z}^4_{\lambda_{\delta},i} \leq A^{1/4}E\bar{X}_k^4$ by (\ref{CD31}) and $x^2 \bar{\mathcal{L}}_{4n} \leq 1/A$. Then 
\begin{eqnarray}\label {69tq}
P\Big (  x B_n^{-2}\Big |\sum_{k \in d}  u_k \Big |> l \Big)  \leq  2e^{-t l} e^{A^{1/2} t^2 A^{1/4}x^2 \bar{\mathcal{L}}_{4n}} \leq  2 e^{- (\log A)l^2/4+l^2}.
\end{eqnarray}
Applying (\ref{69tq}) to (\ref{68tqa}), we have
\begin{eqnarray}\label {72tq}
K_2 \leq  \sum_{l =2}^{[2x]} e^{-l^2} \leq A. 
\end{eqnarray}
For $K_3$, if $m x <  x B_n^{-2}\sum_{k \in d} u_k \leq (m+1)x$, then $m   -2/x^2 <  B_n^{-2}\sum_{k=1}^n u_k \leq m +1 +2/x^2$. By the assumptions of the lemma, $\theta \leq1/(m-1)$ for $m x <  x B_n^{-2}\sum_{k \in d} u_k \leq (m+1)x$ with $m \geq 2$. Hence
\begin{eqnarray}
K_3 \leq \sum_{m=2}^{\infty} e^{10 x^2  (m+1)^2 /(m-1)  } P\Big (   B_n^{-2}\sum_{k \in d}  u_k > m \Big) . \nonumber
\end{eqnarray}
Let $t = (\log A^{1/2})x^2$. Then
\begin{eqnarray}
P\Big (   B_n^{-2}\sum_{k \in d} u_k  > m \Big)  &\leq&   e^{-tm}\prod_{k \in d} \Big (1+\frac{1}{2} e^{t B_n^{-2} |u_k|}t^2 B_n^{-4}E u_k^2\Big)\nonumber\\
&\leq&  e^{-(\log A)m x^2/2+x^2}  . \nonumber
\end{eqnarray}
Hence
\begin{eqnarray}\label {75tq}
K_3 \leq A \sum_{m=2}^{\infty} e^{- m x^2} \leq A. 
\end{eqnarray}
Therefore, (\ref{66tq}) follows from (\ref{67tq}), (\ref{68tq}), (\ref{72tq}) and (\ref{75tq}). 

Now we estimate the first term on the right-hand side of (\ref{44th}). Since $\theta \leq 1/2$, then for each $i$ and $d = \{ 1,...,n \}\backslash \{  i \}$, 
\begin{eqnarray}\label {54th}
E u_i^2  e^{5\theta \bigtriangledown_{n}}  &\leq&   E  u_i^2 e^{10   x^2 B_n^{-4} u_i^2 +10  x^2 B_n^{-4} (\sum_{k \in d} u_k)^2} \nonumber\\
&=& Ee^{10   x^2 B_n^{-4} (\sum_{k \in d} u_k)^2}E u_i^2 e^{10   x^2 B_n^{-4} u_i^2}. 
\end{eqnarray}
Since $x^2 B_n^{-2}|u_i| \leq 1$, then applying (\ref{66tq}) to (\ref{54th}), we have
\begin{eqnarray}\label {55th}
\frac{x^2}{ B_n^{4}}\sum_{i=1}^n E u_i^2  e^{5\theta \bigtriangledown_{n}}  \leq \frac{A x^2}{B_n^4}\sum_{i=1}^n E u_i^2 \leq A x^2 \bar{\mathcal{L}}_{4n}. 
\end{eqnarray}

Next we estimate the second term on the right-hand side of (\ref{44th}). For each $i$ and $j$ with $d = \{ 1,...,n \}\backslash \{  i,j\}$, 
\begin{eqnarray}
&& E u_i u_j e^{5\theta \bigtriangledown_{n}} \nonumber\\
&=& E u_i u_j e^{5 \theta x^2 B_n^{-4}(u_i+u_j)^2 + 10 \theta x^2 B_n^{-4}(u_i+u_j)\sum_{k \in d}u_k+ 5 \theta x^2 B_n^{-4} (\sum_{k \in d} u_k   )^2}. \nonumber
\end{eqnarray}
Since $|e^s-1-s| \leq e^{0 \vee s}s^2/2$ for $s =10 \theta x^2 B_n^{-4}(u_i+u_j)\sum_{k \in d}u_k $ and since $10 \theta x^2 B_n^{-4}|(u_i+u_j)\sum_{k \in d}u_k| \leq 10B_n^{-2}|\sum_{k \in d}u_k| $, then
\begin{eqnarray}\label {58th}
 \frac{x^2}{ B_n^{4}}\sum_{i \neq j} E u_i u_j e^{5\theta \bigtriangledown_{n}} &\leq& \frac{x^2}{ B_n^{4}}\sum_{i \neq j} E  u_i u_j e^{5 \theta x^2 B_n^{-4}(u_i+u_j)^2 } e^{5 \theta x^2 B_n^{-4} (\sum_{k \in d} u_k   )^2} \nonumber\\
&& +\frac{10 \theta x^4}{ B_n^{8}}\sum_{i \neq j}E \bigg  \{  u_i u_j (u_i+u_j)e^{5 \theta x^2 B_n^{-4}(u_i+u_j)^2 }\nonumber\\
&& \times   \Big (\sum_{k \in d}u_k \Big )  e^{  5 \theta x^2 B_n^{-4}(\sum_{k \in d}u_k)^2}\bigg\}\nonumber\\
&& +\frac{50 \theta^2 x^6}{ B_n^{12}}\sum_{i \neq j}  E  |u_i u_j| (u_i+u_j)^2 e^{5   x^2 B_n^{-4}(u_i+u_j)^2 } \nonumber\\
&& \times   E  \Big (\sum_{k \in d}u_k \Big )^2  e^{ 10 B_n^{-2}  |\sum_{k \in d}u_k|+5  x^2 B_n^{-4}(\sum_{k \in d}u_k)^2}\nonumber\\
&:=& L_1 +L_2 +L_3. 
\end{eqnarray}
For $L_1$, since $|e^s-1| \leq e^{s \vee 0}|s|$, $E u_i u_j=0$ and $\theta \leq 1/2$, then 
\begin{eqnarray}
\Big |E \Big (u_i u_j e^{5 \theta x^2 B_n^{-4}(u_i+u_j)^2 }\Big | u_k , k \in  d \Big )\Big | \leq A  x^2 B_n^{-4} E |u_i u_j| (u_i+u_j)^2.  \nonumber
\end{eqnarray}
Then by (\ref{66tq}),
\begin{eqnarray} \label {59th}
L_1 &=&  \frac{x^2}{ B_n^{4}}\sum_{i \neq j}  E \Big \{e^{5 \theta x^2 B_n^{-4} (\sum_{k \in d} u_k   )^2} E \Big (u_i u_j e^{5 \theta x^2 B_n^{-4}(u_i+u_j)^2 }\Big | u_k , k \in  d \Big ) \Big \}  \nonumber\\
&\leq&  \frac{A x^4 }{B_n^{8}} \sum_{i \neq j} E |u_iu_j| (u_i +u_j)^2 \leq A x^2 \bar{\mathcal{L}}_{4n}. 
\end{eqnarray}
For $L_2$, since $|e^s-1 |\leq e^{0 \vee s}|s|$ and $Eu_i u_j(u_i+u_j)=0$,  then
\begin{eqnarray}\label {83tqa}
\Big |E \Big (u_i u_j (u_i+u_j)e^{5 \theta x^2 B_n^{-4}(u_i+u_j)^2 }\Big | u_k, k \in d \Big )\Big | \leq A x^2 B_n^{-4}   E |u_i u_j| |u_i +u_j|^3.
\end{eqnarray}
By H\"{o}lder's inequality and (\ref{66tq}),
\begin{eqnarray}\label {59tk}
&& \frac{x}{B_n^2} E \Big |  \sum_{k \in d}u_k \Big | e^{   5 \theta x^2 B_n^{-4}(\sum_{k \in d}u_k)^2}\nonumber\\
&\leq& \bigg \{\frac{x^2}{B_n^4} E \Big (\sum_{k \in d}u_k \Big)^2 \bigg \}^{1/2} \Big \{E ^{10 \theta x^2 B_n^{-4}(\sum_{k \in d}u_k)^2}\Big \}^{1/2} \leq A. 
\end{eqnarray}
By (\ref{83tqa}) and (\ref{59tk}), 
\begin{eqnarray}\label {60tk}
L_2  &=& \frac{10 \theta x^4}{ B_n^{8}}\sum_{i \neq j}   E \bigg \{ \Big (\sum_{k \in d}u_k \Big )  e^{  5 \theta x^2 B_n^{-4}(\sum_{k \in d}u_k)^2}\nonumber\\
&& \times   E \Big (u_i u_j (u_i+u_j)e^{5 \theta x^2 B_n^{-4}(u_i+u_j)^2 }\Big | u_k, k \in d \Big ) \bigg \}\nonumber\\   \nonumber\\
&\leq&   \frac{ A x^5}{ B_n^{10}}\sum_{i \neq j} E |u_i u_j| |u_i +u_j|^3 \leq A x^2 \bar{\mathcal{L}}_{4n}. 
\end{eqnarray}
For $L_3$, 
\begin{eqnarray}\label {87tqa}
 E \Big ( |u_i u_j| (u_i+u_j)^2 e^{5 \theta x^2 B_n^{-4}(u_i+u_j)^2 }\Big | u_k , k \in d \Big ) \leq A E  ( |u_i u_j| (u_i^2+ u_i^2 )). 
\end{eqnarray}
Applying H\"{o}lder's inequality twice and by (\ref{66tq}),
\begin{eqnarray}  
&& \frac{x^2}{B_n^4} E \Big ( \sum_{k \in d}u_k \Big )^2 e^{20 B_n^{-2}  |\sum_{k \in d}u_k|+5 \theta x^2 B_n^{-4}(\sum_{k \in d}u_k)^2}\nonumber\\
&\leq& \frac{x^2}{B_n^4}\bigg \{ E \Big (\sum_{k \in d}u_k \Big)^4 e^{ 40 B_n^{-2}  |\sum_{k \in d}u_k|}\bigg \}^{1/2} \Big \{E ^{10 \theta x^2 B_n^{-4}(\sum_{k \in d}u_k)^2}\Big \}^{1/2}\nonumber\\
&\leq& A \bigg \{ \frac{x^8}{B_n^{16}} E \Big (\sum_{k \in d}u_k \Big)^8 \bigg \}^{1/4} \Big  \{E e^{80 B_n^{-2}  |\sum_{k \in d}u_k|}\Big  \}^{1/4} . \nonumber
\end{eqnarray}
Since $|e^s-1-s|\leq e^{0 \vee s}s^2/2$ and $E u_k=0$, then
\begin{eqnarray} 
E e^{ 80 B_n^{-2}|\sum_{k \in d}u_k|} &\leq&  E e^{ 80 B_n^{-2}\sum_{k \in d}u_k} + E e^{- 80 B_n^{-2}\sum_{k \in d}u_k} \nonumber\\
&\leq& 2 \prod_{k \in d} (1+A  B_n^{-4} E u_k^2) \leq 2 e^{A B_n^{-4}\bar{\mathcal{L}}_{4n}} \leq 2 e. \nonumber
\end{eqnarray}
It is easy to check that $x^8 B_n^{-16} E  (\sum_{k \in d}u_k  )^8 \leq A x^2 \bar{\mathcal{L}}_{4n} \leq 1$. Then from the above, we have
\begin{eqnarray}\label {90tqa}
 \frac{x^2}{B_n^4} E \Big ( \sum_{k \in d}u_k \Big )^2 e^{20  B_n^{-2}  |\sum_{k \in d}u_k|+5 \theta x^2 B_n^{-4}(\sum_{k \in d}u_k)^2} \leq A. 
\end{eqnarray}
By (\ref{87tqa}) and (\ref{90tqa}), 
\begin{eqnarray}\label {63th}
L_3\leq \frac{Ax^4}{B_n^8} \sum_{ i\neq j} E (|u_i u_j|( u_j^2 +   u_i^2))   \leq A x^2 \bar{\mathcal{L}}_{4n}. 
\end{eqnarray}
Combining (\ref{58th}), (\ref{59th}), (\ref{60tk}) and (\ref{63th}), we obtain
\begin{eqnarray}\label {64th}
\frac{x^2}{ B_n^{4}}\sum_{i \neq j} E u_i u_j e^{5\theta \bigtriangledown_{n}} \leq A x^2 \bar{\mathcal{L}}_{4n}. 
\end{eqnarray}
Therefore, the lemma follows from (\ref{44th}), (\ref{55th}) and (\ref{64th}).
\end{proof}


\end{document}